\newif\ifconf
\conffalse
\conftrue %last one wins

\documentclass[10pt]{article}
\usepackage{amssymb}
\usepackage{amsthm}
\usepackage{amsmath}
\usepackage{fullpage}
\usepackage{epsfig}
\usepackage{verbatim}
\usepackage{hyperref,txfonts}

\DeclareMathOperator{\E}{\mathbb{E}}

\newcommand{\M}{\mathcal{M}}

\newcommand{\A}{\mathcal{A}}

\newcommand{\R}{\mathbb{R}}
\newcommand{\Z}{\mathbb{Z}}

\newcommand{\op}{\mathcal{E}}

\newcommand{\Le}{\Biggl}
\newcommand{\Ri}{\Biggr}

\newcommand{\e}{\varepsilon}

\theoremstyle{plain}
\newtheorem{theorem}{Theorem}[section]
\newtheorem{lemma}[theorem]{Lemma}

\theoremstyle{definition}

\newtheorem{definition}[theorem]{Definition}

\newtheorem{remark}[theorem]{Remark}

\date{}
\begin{document}

\title{Scaled Enflo type is equivalent to Rademacher type}

\author{  Manor Mendel\\ California Institute of Technology
    \and Assaf Naor\\ Microsoft Research
} \maketitle

\begin{abstract}
We introduce the notion of {\em scaled Enflo type} of a metric
space, and show that for Banach spaces, scaled Enflo type $p$ is
equivalent to Rademacher type $p$.
\end{abstract}

\section{Introduction}

Recall that a Banach space $X$ is said to have Rademacher type
$p>0$ (see~\cite{Maurey03}) if there exists a constant $T<\infty$
such that for every $x_1,\ldots,x_n\in X$,
\begin{eqnarray}\label{eq:def type}
\E_\e\Le\|\sum_{j=1}^n \e_j x_j\Ri\|_X^p\le T^p\sum_{j=1}^n
\|x_j\|_X^p,
\end{eqnarray}
where here, and in what follows, $\E_\e$ denotes the expectation
with respect to uniformly chosen $\e=(\e_1,\ldots,\e_n)\in
\{-1,1\}^n$. The infimum over all constants $T$ for
which~\eqref{eq:def type} holds is denoted $T_p(X)$.

 Motivated by the search for concrete
versions of Ribe's theorem~\cite{Ribe76} for various fundamental
local properties of Banach spaces (see the discussion
in~\cite{Bourgain86-trees,NPSS04,MN05-cotype}), several
researchers proposed non-linear notions of type, which make sense
in the setting arbitrary metric spaces
(see~\cite{Enflo78,BMW86,Ball92}). In particular, following
Enflo~\cite{Enflo78} we say that a metric space $(\M,d_\M)$ has
{\em Enflo type $p$} if there exists a constant $K$ such that for
every $n\in \mathbb N$ and  every $f:\{-1,1\}^n\to \M$,
\begin{eqnarray}\label{eq:enflo type}
\E_\e d_\M(f(\e),f(-\e))^p\le T^p\sum_{j=1}^n \E_\e
d_\M\left(f(\e_1,\ldots,\e_{j-1},\e_j,\e_{j+1},\ldots,\e_n),f(\e_1,\ldots,\e_{j-1},-\e_j,\e_{j+1},\ldots,\e_n)\right)^p.
\end{eqnarray}

For Banach spaces~\eqref{eq:def type} follows from~\eqref{eq:enflo
type} by considering the function $\e\mapsto \sum_{j=1}^n \e_j
x_j$. The question whether in the category of Banach spaces
Rademacher type $p$ implies Enflo type $p$ was posed by Enflo
in~\cite{Enflo78}, and in full generality remains open.
In~\cite{Pisier86} Pisier showed that if a Banach space has
Rademacher $p$ then it has Enflo type $p'$ for every $p'<p$ (see
also the work of Bourgain, Milman and Wolfson~\cite{BMW86} for a
similar result which holds for a another notion of non-linear
type). In~\cite{NS02} it was shown that for UMD Banach spaces
(see~\cite{Burk01}) Rademacher type $p$ is equivalent to Enflo
type $p$.

Motivated by our recent work on metric cotype~\cite{MN05-cotype},
we introduce below the notion of {\em scaled Enflo type} of a
metric space (which is, in a sense, ``opposite" to the notion of
metric cotype defined in~\cite{MN05-cotype}), and show that for
Banach spaces, scaled Enflo type $p$ is equivalent to Rademacher
type $p$. This settles the long standing problem of finding a
purely metric formulation of the notion of type (though Enflo's
problem described above remains open). Modulo some of the results
of~\cite{MN05-cotype}, the proof of our main theorem is very
simple.

\begin{definition}[Scaled Enflo type]\label{del:scaled}
Let $(\M,d_\M)$ be a metric space and $p>0$. We say that $\M$ has
{\em scaled Enflo type} $p$ with constant $\tau$ if for every
integer $n$ there exists an even integer $m$ such that for every
$f:\Z_m^n\to \M$,
\begin{eqnarray}\label{eq:def scaled}
\E_\e\int_{\Z_m^n}d_\M\left(f\left(x+\frac{m}{2}\e\right),f(x)\right)^pd\mu(x)\le
\tau^pm^p\sum_{j=1}^n\int_{\Z_m^n}d_\M\left(f(x+e_j),f(x)\right)^pd\mu(x),
\end{eqnarray}
where $\mu$ is the uniform probability measure on $\Z_m^n$, and
$\{e_j\}_{j=1}^n$ is the standard basis of $\R^n$. The infimum
over all constants $\tau$ for which~\eqref{eq:def scaled} holds is
denoted $\tau_p(\M)$.
\end{definition}

\begin{theorem}\label{thm:main} Let $X$ be a Banach space and $p\in [1,2]$. Then $X$ has Rademacher type $p$ if and only if $X$ has
scaled Enflo type $p$. More precisely,
$$
\frac{1}{2\pi}T_p(X)\le \tau_p(X)\le 5T_p(X).
$$
\end{theorem}

\section{Proof of Theorem~\ref{thm:main}}

We start by showing that scaled Enflo type $p$ implies Rademacher
type $p$.

\begin{lemma}\label{lem:Enflo implies Rademacher}
Let $X$ be a Banach space and $p\in [1,2]$. Then $ T_p(X)\le
2\pi\tau_p(X)$.
\end{lemma}

\begin{proof}
Let $X$ be a Banach space and assume that $\tau_p(X)<\infty$ for
some $p\in [1,2]$. Fix $\tau>\tau_p(X)$, $v_1,\ldots,v_n\in X$,
and let $m$ be an even integer. Define $f:\Z_m^n\to X$ by $
f(x_1,\ldots,x_n)=\sum_{j=1}^n e^{\frac{2\pi i x_j}{m}}v_j$. Then
\begin{eqnarray}\label{eq:reverse cotype}
\sum_{j=1}^n\int_{\Z_m^n}\left\|f\left(x+e_j\right)-f(x)\right\|_X^pd\mu(x)=\left|e^{\frac{2\pi
i}{m}}-1\right|^p\cdot\sum_{j=1}^n \|v_j\|_X^p\le
\left(\frac{2\pi}{m}\right)^p\cdot\sum_{j=1}^n \|v_j\|_X^p,
\end{eqnarray}
and
\begin{eqnarray}\label{eq:reverse cotype2}
\E_\e\int_{\Z_m^n}\Le\|f\left(x+\frac{m}{2}\e\right)-f(x)\Ri\|_X^pd\mu(x)=
2^p\int_{\Z_m^n}\Le\|\sum_{j=1}^n e^{\frac{2\pi i x_j}{m}}v_j
\Ri\|_X^pd\mu(x).
\end{eqnarray}
We recall the {\em contraction principle} (see~\cite{LT91}), which
states that for every $a_1,\ldots, a_n\in \R$,
$$
\E_\e\Le\|\sum_{j=1}^n \e_ja_jv_j\Ri\|_X^p\le \left(\max_{1\le
j\le n} |a_j|\right)^p\cdot \E_\e\Le\|\sum_{j=1}^n
\e_jv_j\Ri\|_X^p.
$$
Thus,
\begin{eqnarray}\label{eq:before plug}
\int_{\Z_m^n}\Le\|\sum_{j=1}^n e^{\frac{2\pi i x_j}{m}}v_j
\Ri\|_X^pd\mu(x)=\int_{\Z_m^n}\E_\e\Le\|\sum_{j=1}^n e^{\frac{2\pi
i }{m}\left(x_j+\frac{m(1-\e_j)}{4}\right)}v_j
\Ri\|_X^pd\mu(x)=\int_{\Z_m^n}\E_\e\Le\|\sum_{j=1}^n
\e_je^{\frac{2\pi i x_j}{m}}v_j \Ri\|_X^pd\mu(x)\ge
\frac{1}{2^p}\E_\e\Le\|\sum_{j=1}^n \e_jv_j \Ri\|_X^p.
\end{eqnarray}
Combining~\eqref{eq:reverse cotype}, \eqref{eq:reverse cotype2}
and~\eqref{eq:before plug} yields the required result.
\end{proof}

Let $X$ be a Banach space with type $p$, $m$ an integer divisible
by $4$, and $k$ an odd integer. Fix $f:\Z_m^n\to X$ and $\e\in
\{-1,1\}^n$. Define $\A^{(k)}f:\Z_m^n\to X$ by
$$
\A^{(k)}f(x)= \frac{1}{k^n} \sum_{z\in (-k,k)^n\cap
(2\Z)^n}f(x+z).
$$

\begin{lemma}\label{lem:approxA} For $p\ge 1$ and every
$f:\Z_m^n\to X$
$$
\int_{\Z_m^n} \left\|\A^{(k)}f(x)-f(x)\right\|_X^pd\mu(x)\le
(k-1)^{p}n^{p-1}\sum_{j=1}^n\int_{\Z_m^n}\|f(x+e_j)-f(x)\|_X^pd\mu(x).
$$
\end{lemma}

\begin{proof} For every $t\in \R$ let $s(t)$ be the sign of $t$
(with convention that $s(0)=0$). For every $z\in \Z_m^n$,
\begin{eqnarray*}
\|f(x+z)-f(x)\|_X^p\le  \|z\|_1^{p-1}\cdot\sum_{j=1}^n
\sum_{\ell=1}^{|z_j|}\Le\|f\Le(x+\sum_{t=1}^{j-1}z_t e_t+\ell\cdot
s(z_j)\cdot e_j\Ri)-f\Le(x+\sum_{t=1}^{j-1}z_t e_t+ (\ell-1)\cdot
s(z_j)\cdot e_j\Ri)\Ri\|_X^p.
\end{eqnarray*}
Observe that since $k$ is odd, $|(-k,k)^n\cap (2\Z)^n|=k^n$. Thus
\begin{eqnarray*}
&&\!\!\!\!\!\!\!\!\!\!\!\!\!\!\!\!\int_{\Z_m^n}
\left\|\A^{(k)}f(x)-f(x)\right\|_X^pd\mu(x)\le\frac{1}{k^n}
\sum_{z\in (-k,k)^n\cap
(2\Z)^n}\int_{\Z_m^n}\|f(x+z)-f(x)\|_X^pd\mu(x)\\
&\le& \frac{1}{k^n}\sum_{z\in (-k,k)^n\cap
(2\Z)^n}\int_{\Z_m^n}\|z\|_1^{p-1}\sum_{j=1}^n
\sum_{\ell=1}^{|z_j|}\Le\|f\Le(x+\sum_{t=1}^{j-1}z_t e_t+\ell
s(z_j)
e_j\Ri)-f\Le(x+\sum_{t=1}^{j-1}z_t e_t+(\ell-1) s(z_j) e_j\Ri)\Ri\|_X^pd\mu(x)\\
&\le& \frac{1}{k^n}\sum_{z\in (-k,k)^n\cap (2\Z)^n}\sum_{j=1}^n
\|z\|_1^{p-1}|z_j|\int_{\Z_m^n}\|f(y+s(z_j)e_j)-f(y)\|_X^pd\mu(x)\\
&\le&
(k-1)^{p}n^{p-1}\sum_{j=1}^n\int_{\Z_m^n}\|f(x+e_j)-f(x)\|_X^pd\mu(x).
\end{eqnarray*}
\end{proof}

\begin{proof}[Proof of theorem~\ref{thm:main}] Fix an odd integer $k\in \mathbb N$, with $k<
\frac{m}{2}$. As in~\cite{MN05-cotype},  given $j\in
\{1,\ldots,n\}$ we define $S(j,k)\subseteq \Z_m^n$ by
$$
S(j,k)\coloneqq\left\{x\in [-k,k]^n\subseteq \Z_m^n:\ y_j\equiv
0\mod 2\ \mathrm{and}\ \forall\  \ell\neq j,\ y_\ell\equiv 1\mod
2\right\}.
$$
For $f:\Z_m^n\to X$ we define
\begin{eqnarray}\label{eq:def Ek}
\op_j^{(k)}f(x)=\left(f*\frac{\mathbf{1}_{S(j,k)}}{\mu(S(j,k))}\right)(x)=\frac{1}{\mu(S(j,k))}\int_{S(j,k)}
f(x+y)d\mu(y).
\end{eqnarray}
In~\cite{MN05-cotype} (see equation (39) there) it is shown that
for every $x\in \Z_m^n$ and $\e\in \{-1,1\}^n$,
$$
\left(\frac{k}{k+1}\right)^{n-1}\left(\A^{(k)}f(x+\e)-\A^{(k)}f(x-\e)\right)=\sum_{j=1}^n
\e_j\left[\op_j^{(k)}f(x+e_j)-\op_j^{(k)}f(x-e_j)\right]+U(x,\e)+V(x,\e),
$$
where, by inequalities (41) and (42) in~\cite{MN05-cotype}, for
every $\e\in \{-1,1\}^n$,
$$
\max\left\{\int_{\Z_m^n}\|U(x,\e)\|_X^pd\mu(x),\int_{\Z_m^n}\|V(x,\e)\|_X^pd\mu(x)\right\}\le
\frac{8^{p}n^{2p-1}}{k^p}\sum_{j=1}^n\int_{\Z_m^n}\left\|f(x+e_j)-f(x)\right\|_X^p.
$$
Thus, for every $T>T_p(X)$,
\begin{eqnarray}\label{eq:step}
&&\!\!\!\!\!\!\!\!\!\!\!\!\!\!\!\!\nonumber\left(\frac{k}{k+1}\right)^{p(n-1)}\E_\e\int_{\Z_m^n}
\left\|\A^{(k)}f(x+\e)-\A^{(k)}f(x-\e)\right\|_X^pd\mu(x)\\\nonumber\label{eq:from
cotype paper}&\le& 3^{p-1}\E_\e\int_{\Z_m^n}\Le\|\sum_{j=1}^n
\e_j\left[\op_j^{(k)}f(x+e_j)-\op_j^{(k)}f(x-e_j)\right]\Ri\|_X^pd\mu(x)+\frac{24^p
n^{2p-1}}{k^p}\sum_{j=1}^n\int_{\Z_m^n}\|f(x+e_j)-f(x)\|_X^pd\mu(x)\\\nonumber
&\le& 3^{p-1}T^p\sum_{j=1}^n
\int_{\Z_m^n}\left\|\op_j^{(k)}f(x+e_j)-\op_j^{(k)}f(x-e_j)\right\|_X^pd\mu(x)+\frac{24^p
n^{2p-1}}{k^p}\sum_{j=1}^n\int_{\Z_m^n}\|f(x+e_j)-f(x)\|_X^pd\mu(x)\\\label{eq:averaging}&\le&
3^{p-1}T^p\sum_{j=1}^n
\int_{\Z_m^n}\left\|f(x+e_j)-f(x-e_j)\right\|_X^pd\mu(x)+\frac{24^p
n^{2p-1}}{k^p}\sum_{j=1}^n\int_{\Z_m^n}\|f(x+e_j)-f(x)\|_X^pd\mu(x)\\
&\le& \left(\frac{6^p}{3}T^p+\frac{24^p
n^{2p-1}}{k^p}\right)\sum_{j=1}^n\int_{\Z_m^n}\|f(x+e_j)-f(x)\|_X^pd\mu(x),
\end{eqnarray}
where in~\eqref{eq:averaging} we used the fact that $\op_j^{(k)}$
is an averaging operator, and hence has norm $1$.

On the other hand
\begin{eqnarray}
&&\!\!\!\!\!\!\!\!\!\!\nonumber\E_\e\int_{\Z_m^n}\left\|f\left(x+\frac{m}{2}\e\right)-f(x)\right\|_X^pd\mu(x)
\le
3^{p-1}\E_\e\int_{\Z_m^n}\left\|\A^{(k)}f\left(x+\frac{m}{2}\e\right)-\A^{(k)}f(x)\right\|_X^pd\mu(x)+\\\nonumber&\phantom{\le}&
3^{p-1}\E_\e\int_{\Z_m^n}\left\|f\left(x+\frac{m}{2}\e\right)-\A^{(k)}f\left(x+\frac{m}{2}\e\right)\right\|_X^pd\mu(x)+3^{p-1}
\E_\e\int_{\Z_m^n}\left\|\A^{(k)}f(x)-f(x)\right\|_X^pd\mu(x)\\\nonumber
&\le&
3^{p-1}\left[\left(\frac{m}{4}\right)^{p-1}\E_\e\int_{\Z_m^n}\sum_{t=1}^{m/4}\left\|\A^{(k)}f\left(x+2t\e\right)-\A^{(k)}f(x+(2t-2)\e)\right\|_X^pd\mu(x)+2
\E_\e\int_{\Z_m^n}\left\|\A^{(k)}f(x)-f(x)\right\|_X^pd\mu(x)\right]
\\\label{eq:the lemma}
&\le&3^{p-1}\left[\left(\frac{m}{4}\right)^{p}\E_\e\int_{\Z_m^n}\sum_{t=1}^{m/4}\left\|\A^{(k)}f\left(y+\e\right)-\A^{(k)}f(x-\e)\right\|_X^pd\mu(x)+2k^pn^{p-1}
\sum_{j=1}^n\int_{\Z_m^n}\|f(x+e_j)-f(x)\|_X^pd\mu(x)\right]\\\label{eq:the
step} &\le&
\left[3^{p-1}\left(\frac{m}{4}\right)^{p}\left(1+\frac{1}{k}\right)^{p(n-1)}\left(\frac{6^p}{3}T^p+\frac{24^p
n^{2p-1}}{k^p}\right)+\frac{2(3kn)^p}{3n}\right]\sum_{j=1}^n\int_{\Z_m^n}\|f(x+e_j)-f(x)\|_X^pd\mu(x)\\\label{eq:big
m} &\le&
5^pm^pT^p\sum_{j=1}^n\int_{\Z_m^n}\|f(x+e_j)-f(x)\|_X^pd\mu(x),
\end{eqnarray}
where in~\eqref{eq:the lemma} we used Lemma~\ref{lem:approxA},
in~\eqref{eq:the step} we used~\eqref{eq:step}, and~\eqref{eq:big
m} is true if  $4n^{2-1/p}\le k\le \frac{3m}{2n^{1-1/p}}$, which
is a valid choice of $k$ if $m\ge 3n^{3-2/p}$.
\end{proof}

\begin{remark} If a metric space has Enflo type $p$ then it also
has scaled Enflo type $p$. This follows from a straightforward
modification of Lemma 2.4 in~\cite{MN05-cotype}. We do not know if
scaled Enflo type $p$ implies Enflo type $p$. In the category of
Banach spaces, a positive answer to this question would show that
Enflo type $p$ is equivalent to Rademacher type $p$, resolving
positively Enflo's problem~\cite{Enflo78}. We do know that for
Banach spaces, scaled Enflo type $p$ implies Enflo type $p'$
for all $p'<p$, and that scaled Enflo type and Enflo type coincide
for UMD Banach spaces.
\end{remark}

\begin{remark}
The idea of scaling by $\frac{m}{2}$ in the definition of scaled
Enflo type originates from the definition of {\em metric cotype}
introduced in~\cite{MN05-cotype}, which involves a similar scaling
procedure. In the case on non-linear type it is possible that this
scaling is not necessary, i.e. that Enflo type is equivalent to
Rademacher type. However, as shown in~\cite{MN05-cotype}, in the
context of metric cotype the scaling {\em is necessary}- we refer
to~\cite{MN05-cotype} for more details.
\end{remark}

%\cite{Enflo78}\cite{Pisier86}\cite{BMW86}\cite{Ribe76}\cite{Bourgain86-trees}\cite{NS02}\cite{Ball92}\cite{MN05-cotype}\cite{NPSS04}\cite{Burk01}

%\bibliographystyle{abbrv}
%\bibliography{cotype}

\def\cprime{$'$}

\end{document}